\RequirePackage{ifpdf}
\ifpdf 
\documentclass[pdftex]{sigma}
\else
\documentclass{sigma}
\fi

\def\b{\mathbb }
\def\phi{\varphi }

\begin{document}

\renewcommand{\thefootnote}{$\star$}

\renewcommand{\PaperNumber}{083}

\FirstPageHeading

\ShortArticleName{A Limit Relation for Dunkl--Bessel Functions of Type A and B}

\ArticleName{A Limit Relation for Dunkl--Bessel Functions\\
 of Type A and B\footnote{This paper is a contribution to the Special
Issue on Dunkl Operators and Related Topics. The full collection
is available at
\href{http://www.emis.de/journals/SIGMA/Dunkl_operators.html}{http://www.emis.de/journals/SIGMA/Dunkl\_{}operators.html}}}

\Author{Margit R\"OSLER~$^\dag$ and Michael VOIT~$^\ddag$}

\AuthorNameForHeading{M. R\"osler and M. Voit}

\Address{$^\dag$~Institut f\"ur Mathematik, TU Clausthal,
Erzstr. 1, D-38678 Clausthal-Zellerfeld, Germany}

\EmailD{\href{mailto:roesler@math.tu-clausthal.de}{roesler@math.tu-clausthal.de}}

\Address{$^\ddag$~Fachbereich Mathematik, TU Dortmund, Vogelpothsweg 87,
          D-44221 Dortmund, Germany}
\EmailD{\href{mailto:michael.voit@math.tu-dortmund.de}{michael.voit@math.tu-dortmund.de}}

\ArticleDates{Received October 21, 2008, in f\/inal form November 26,
2008; Published online December 03, 2008}

\Abstract{We prove a limit relation for the Dunkl--Bessel function of type $B_N$
with multiplicity parameters  $k_1$  on the roots $\pm e_i$ and $k_2$ on  $\pm e_i\pm e_j$ where $k_1$ tends to inf\/inity and the arguments
are suitably scaled. It gives a good approximation
in terms of the Dunkl-type Bessel function of type $A_{N-1}$ with multiplicity $k_2.$ For certain values of
$k_2$ an improved estimate is obtained from a corresponding limit relation
for Bessel functions on matrix cones.}


\Classification{33C67; 43A85; 20F55}

\section{Introduction and results}

The power series expansion
\[j_\alpha(z)=\Gamma(\alpha+1)\sum_{n=0}^\infty
\frac{(-1)^n (z/2)^{2n }}{n! \> \Gamma(n+\alpha+1)}\]
of the normalized spherical Bessel functions $j_\alpha(z) = {} _0F_1(\alpha +1; -z^2/4)$ leads immediately to the well
known limit relation
\begin{equation*}
\lim_{\alpha\to\infty} j_\alpha(\sqrt\alpha \cdot z)= e^{-z^2/4}
\qquad \text{for} \quad z\in\b C.
\end{equation*}
For real arguments, this limit relation can be improved  as follows: There exists a
constant $C>0$ such that
\begin{equation}\label{bessel-estimate-1-dim}
\big|j_{\mu -1}(\sqrt\mu \cdot x) - e^{-x^2/4}\big|\le \frac{C}{\mu}\cdot
 \min\big(x^4,1\big)
\qquad\text{for all}\quad x\in\b R, \  \mu >2.
\end{equation}
This is the rank-one specialization of a more general result for Bessel functions on cones of
 positive semidef\/inite matrices obtained in  \cite[Theorem~3.6]{RV2}.   For related asymptotic results for
one-variable Bessel functions $j_\alpha$ as $\alpha\to\infty$ we refer to~\cite{W}.
 In~\cite{V}, a  variant of estimate~(\ref{bessel-estimate-1-dim}) was used to derive a law of large
numbers for radial random walks on $\b R^p$ where the time parameter of the
walks as well as the dimension $p$ tend to inf\/inity. In this case, Bessel functions of index
 $\alpha=\frac{p}{2}-1$ come in as the radial
parts of the complex exponential functions, and estimate~(\ref{bessel-estimate-1-dim}) leads
to limit theorems for radial random walks. This approach was extended in
\cite{RV2} to radial random walks on matrix spaces over one of the (skew-) f\/ields
$\mathbb F=\mathbb R, \mathbb C, \mathbb H.$

In the present note we  derive a further multidimensional extension of
(\ref{bessel-estimate-1-dim}),
namely for certain Bessel functions of Dunkl-type; see \cite{DX,R1} for an introduction
 to Dunkl theory and~\cite{O} for  the associated  Bessel functions.
More precisely, we shall prove that under suitable normalization of the arguments, the  Bessel
 function  $J_{(k_1,k_2)}^B$
of type $B_q$ converges  to the Bessel function $J_{k_2}^A$
of type $A_{q-1}$ where the multiplicity parameter $k_2$ (on the roots $\pm e_i\pm e_j$) is f\/ixed
and $k_1$ (on the roots $\pm e_i$) tends to inf\/inity. The obtained
estimate is optimal for small arguments, but for large arguments it is weaker than~(\ref{bessel-estimate-1-dim}) and
the corresponding result in~\cite{RV2} for Bessel functions of matrix argument. This is due to the fact that the proofs
of  (\ref{bessel-estimate-1-dim}) and its matrix version in~\cite{RV2}
rely on some explicit integral representation of the Bessel functions which is
as far not available for Dunkl-type Bessel functions in general.  Nevertheless, our
limit result should be of some interest in its own.
In the following, we state our limit results. Proofs will be
given in Sections~\ref{sec2} and~\ref{sec3}.

We f\/irst recapitulate the necessary facts from Dunkl theory. We
shall not go into details but refer the  reader to \cite{DX,R1} and \cite{RV2} for more background.  For multivariable hypergeometric functions, see
e.g.~\cite{BF1,GR} and \cite{Ka}.
For a reduced root system
 $R\subset \b R^N$ and a multiplicity
 function $k:R\to [0,\infty)$  (i.e.~$k$ is  invariant under the action of the corresponding ref\/lection group),
we denote by $E_k$ the corresponding Dunkl kernel and by
$J_k$ the Bessel function associated with $R$ and $k$ which is given by
\[ J_k(x,y) = \frac{1}{|W|} \sum_{w\in W} E_k(wx,y),\]
where the sum is over the underlying ref\/lection group $W.$
Bessel functions associated with root systems generalize
 the spherical functions of f\/lat symmetric spaces which occur for
 crystallographic root systems
 and specif\/ic discrete values of $k$.  We shall be concerned with Bessel functions associated
with root systems of type $A$ and $B$ which can be expressed in terms of Jack polynomial series.
 To be precise, let  $C_\lambda^\alpha$ denote
the Jack polynomials of index $\alpha >0$ which are indexed by partitions
$\lambda = (\lambda_1\geq  \cdots \geq \lambda_N)\in \mathbb N_0^N,$  see \cite{Sta}. The $C_\lambda^\alpha$
 are homogeneous of degree  $|\lambda|= \lambda_1 + \cdots + \lambda_N$ and can  be normalized such that
\begin{equation}\label{jacksum} (x_1 + \cdots + x_N)^k = \sum_{|\lambda|=k} C_\lambda^\alpha(x)
 \qquad \text{for all}\quad k \in \mathbb N_0;\end{equation}
this normalization will be adopted here.
For root system $A_{N-1}= \{\pm(e_i-e_j):\, i<j\}\subset \b R^N$, the multiplicity $k$ is
a single real parameter. Due to relations~(3.22) and~(3.37)  of \cite{BF2},  the associated Bessel
function can be expressed as
 a generalized $_0F_0$-hypergeometric function,
\begin{equation}\label{power-j-a}
 J_k^A(x,y) ={}_0F_0^\alpha (x,y) := \sum_{\lambda \geq 0} \frac{1}{|\lambda|!}\cdot
\frac{C_\lambda^\alpha(x) C_\lambda^\alpha(y)}{C_\lambda^\alpha({\bf
    1})}\qquad \text{with}\quad {\bf 1}
= (1,\dots, 1),\  \ \alpha = 1/k,
\end{equation}
where  $\lambda \geq 0$ denotes that the sum is taken over all partitions $\lambda=(\lambda_1 \geq \cdots \geq \lambda_N).$
For  root system   $B_N = \{\pm e_i , \,\pm e_i \pm e_j: \, i <j \}$, the
    multiplicity is of the form
 $k=(k_1, k_2)$ where~$k_1$ and~$k_2$ are the values on the roots $\pm e_i$
    and  $\pm e_i \pm e_j$ respectively.
 The associated Bessel function is given by
\begin{equation}\label{power-j-b}
 J_k^B(x,y) = {} _0F_1^\alpha\left(\mu; \frac{x^2}{2}, \frac{y^2}{2}\right)
\qquad \text{with}\quad \alpha = \frac{1}{k_2}, \ \  \mu = k_1 +(N-1)k_2 +\frac{1}{2},
\end{equation}
where $ x^2:= (x_1^2, \ldots, x_N^2)$ and
\[_0F_1^\alpha(\mu; x, y) :=   \sum_{\lambda\geq 0}
 \frac{1}{(\mu)_\lambda^\alpha |\lambda|!}\cdot \frac{C_\lambda^\alpha(x) C_\lambda^\alpha(y)}{C_\lambda^\alpha(\bf{1})}.\]

We denote by $\langle \cdot,\cdot\rangle$ and $\vert \cdot \vert$ the
usual Euclidean scalar product and norm on $\mathbb R^N$.
The main results of this note are as follows:

\begin{proposition}\label{proposition11}
Let $N\geq 2$ and $k_2\geq 0$. Then
there exists a constant $C=C(N,k_2)>0$ such that for all $k_1\ge k_2(N-1)$,
$x,y\in\b R^N$, and $\mu= k_1 + k_2(N-1)+\frac{1}{2}$,
\[ \big\vert J_{(k_1,k_2)}^B(2\sqrt{\mu} x, iy)  -  J_{k_2}^A\big({-}x^2,y^2\big) \big\vert
  \leq  \frac{C}{\mu} \cdot |x|^4|y|^4 \cdot e^{|x|^2|y|^2 }.\]
 \end{proposition}

\noindent
For certain values of $k_2$, this estimate can be improved as follows:

\begin{proposition}\label{proposition12}
Let $N\ge 2$ and $k_2\in\{0, \frac{1}{2}, 1, 2\}$. Then
there exists a constant $C=C(N,k_2)>0$ such that for all $k_1\ge k_2(N-1)$,
$x,y\in\b R^N$, and $\mu$ as above,
\[ \big\vert J_{(k_1,k_2)}^B(2\sqrt{\mu}  x, iy) - J_{k_2}^A\big({-}x^2,y^2\big) \big\vert
  \leq  \frac{C}{\mu} \cdot \min\big(|x|^4|y|^4, 1\big) .\]
 \end{proposition}

In contrast to the previous estimate which is only locally uniform, this estimate is uniform in $x$ and~$y$. For $x$, $y$ close to $0$, it gives the same rate of convergence for $\mu\to\infty$.
We conjecture that
Proposition~\ref{proposition12} is actually correct for all $k_2\geq 0$, and that a similar estimate is valid for the
associated Dunkl kernels of type $A$ and type $B$.

\section{Proof of Proposition \ref{proposition12}}\label{sec2}

\subsection[The case $k_2=0$]{The case $\boldsymbol{k_2=0}$}

The argumentation in this case is dif\/ferent from that in the remaining
cases and based on a~reduction to the rank one case. Indeed,
the Dunkl operators of
type B with multiplicity $(k_1, 0)$ may be regarded as  Dunkl operators
for the ref\/lection group $\mathbb Z_2^N$ and multiplicity $k_1=:k$. Thus the Dunkl kernel $E_{(k,0)}^B$ factorizes as
$ E_{(k,0)}^B(x,y) = \prod_{l=1}^N E_k^{\mathbb Z_2}(x_l, y_l) $ and the associated Bessel function is given by
\[ J_{(k,0)}^B(x,y) = \frac{1}{N!} \sum_{w\in S_N} G_k(wx,y),\qquad \text{where} \qquad G_k(x,y) = \prod_{l=1}^N J_k^{\mathbb Z_2}(x_l,y_l) \]
and $J_k^{\mathbb Z_2}$ is the Bessel function for root system $\mathbb Z_2$ on $\mathbb R$, that is $  J_k^{\mathbb Z_2}(x,y) = j_{k-\frac{1}{2}}(ixy).$
On the other hand, the type $A$ Bessel function with multiplicity $0$ is just
\[ J_0^A(x,y)= \frac{1}{N!} \sum_{w\in S_N} e^{\langle wx,y\rangle}.\]
Thus
\[ \big\vert J_{(k,0)}^B(2\sqrt\mu x, iy) - J_0^A\big({-}x^2,y^2\big)\big\vert \leq \frac{1}{N!}  \sum_{w\in S_N} \big\vert G_k(2\sqrt\mu wx,iy) -e^{-\langle (wx)^2,y^2\rangle}\big\vert.\]
Further, for $x,y\in \mathbb R^N$,
\begin{gather*}
\big\vert G_k(2\sqrt\mu x,iy) - e^{-\langle x^2,y^2\rangle}\big\vert
 = \big\vert \prod_{l=1}^N j_{k-\frac{1}{2}}(2\sqrt\mu  x_ly_l) - \prod_{l=1}^N e^{-x_l^2y_l^2}\big\vert\\
  \phantom{\big\vert G_k(2\sqrt\mu x,iy) - e^{-\langle x^2,y^2\rangle}\big\vert=}{}
  \leq \sum_{l=1}^N \big\vert j_{k-\frac{1}{2}}(2\sqrt\mu  x_ly_l) -e^{-x_l^2y_l^2}\big\vert,
 \end{gather*}
where $\mu = k+ \frac{1}{2}$ and the last inequality is obtained by a telescope argument and the
fact that the factors in both products are bounded by~$1$. By~\eqref{bessel-estimate-1-dim}, the last sum can be estimated by $ \frac{C^\prime}{\mu}\cdot\min(1, |x|^4|y|^4),$ and this  yields the stated result.

\subsection[The cases $k_2=\frac{1}{2}, 1,2$]{The cases $\boldsymbol{k_2=\frac{1}{2}, 1,2}$}

In these cases, the Bessel functions of type B are closely related with Bessel functions on
the matrix cones $\Pi_N=\Pi_N(\b F)$ of positive semidef\/inite $N\times N$ matrices over $\b F=\b R,\b C, \b H$ as explained
in~\cite{R2}. Using this connection, we shall derive Proposition~\ref{proposition12} from a corresponding result for  Bessel
functions on matrix cones in~\cite{RV2}.

We f\/irst recapitulate some facts about Bessel functions of matrix argument. Fix one of the skew-f\/ields
$\b F=\b R,\b C, \b H$ with real dimension $d=1,2,4,$ respectively.
The Bessel functions associated with the cone $\Pi_N = \Pi_N(\mathbb F) $ are def\/ined in terms of its spherical polynomials which are indexed by partitions $\lambda = (\lambda_1\geq \cdots \geq \lambda_N)\in \mathbb N_0^N$ and given by
\[
\Phi_\lambda (X) = \int_{U_N} \Delta_\lambda\big(UXU^{-1}\big)dU,
\]
where $dU$ is the normalized Haar measure of the unitary group $U_N= U_N(\mathbb F)$ and $\Delta_\lambda$ denotes the power function
\[ \Delta_\lambda(X) := \Delta_1(X)^{\lambda_1-\lambda_2}
\Delta_2(X)^{\lambda_2-\lambda_3} \cdot\cdots\cdot
\Delta_N(X)^{\lambda_N}\]
on the vector space $H_N= \{X\in M_N(\mathbb F): X=X^*\}$ of Hermitian  $N\times N$ matrices over $\b F$.
The $\Delta_i(X)$ are the principal minors of the determinant $\Delta(X)$, see
\cite{FK} for details.
There is a~renormalization
$Z_\lambda = c_\lambda \Phi_\lambda$
with constants $c_\lambda >0$ depending on $\Pi_N$
such that
\[
({\rm tr} \,X)^k  =  \sum_{|\lambda|=k} Z_\lambda(X)
\qquad \text{for all}\quad
 k\in \mathbb N_0,
\]
see Section XI.5 of~\cite {FK}. By construction, the $Z_\lambda$ are invariant under
conjugation by $U_N$ and thus
 depend only on the eigenvalues of their argument.
More precisely, for $X\in H_N$ with eigenvalues $x = (x_1, \ldots, x_N)\in
\b R^N$,
\begin{equation*}
Z_\lambda(X) = C_\lambda^\alpha(x) \qquad \text{with}\quad \alpha = \frac{2}{d},
\end{equation*}
where the $C_\lambda^\alpha$ are the Jack polynomials of index
$\alpha$ (cf.~\cite{FK,Ka,R3}).

The Bessel functions on the cone $\Pi_N$ are def\/ined
as  $_0F_1$-hyper\-geometric series in terms of the~$Z_\lambda$, namely
\begin{equation}\label{power-j}
 J_\mu(X) =
 \sum_{\lambda\geq 0} \frac{(-1)^{|\lambda|}}{(\mu)_\lambda|\lambda|!} Z_\lambda(X),
\end{equation}
where the sum is over all partitions $\lambda = (\lambda_1\geq \cdots \geq \lambda_N)\in \mathbb N_0^N$  and $(\mu)_\lambda$ denotes
 the generalized Pochhammer symbol
\[ (\mu)_\lambda =  (\mu)_\lambda ^{2/d} \qquad \text{where}\qquad
(\mu)_\lambda^\alpha :=
 \prod_{j=1}^N \left(\mu-\frac{1}{\alpha}(j-1)\right)_{\lambda_j} \quad (\alpha >0). \]
In \eqref{power-j}, the index $\mu\in \b C$ is supposed to satisfy $(\mu)_\lambda^\alpha \not= 0$ for all $\lambda \geq 0.$ If $N=1,$ then $\Pi_1=\b R_+$  and the Bessel function $\mathcal J_\mu$ is independent of $d$ and given by a usual one-variable Bessel function,
\[ J_\mu\left(\frac{x^2}{4}\right) = j_{\mu-1}(x). \]

There exist commutative convolution algebras (so-called hypergroup structures) on
the co\-ne~$\Pi_N$ with convolutions which depend on the parameter $\mu$ and which have
 the Bessel functions $\phi_Y(X) = \mathcal J_\mu\bigl(\frac{1}{4}YX^2Y\bigr)$, $Y\in \Pi_N,$
 as characters. For details we refer to~\cite{R2}. Moreover, the unitary group $U_N$ acts by the usual conjugation
  $X\mapsto UXU^{-1}$ on $\Pi_N$ as a compact group of hypergroup automorphisms, i.e.,
 these conjugations preserve these convolution structures. As shown in~\cite{R2}, this observation induces a further commutative hypergroup structure on the associatd
orbit space $\Pi_N^{U_N}$ where this space may obviously be identif\/ied with the
 the set
of all possible eigenvalues of matrices from $\Pi_N$ ordered by size, i.e.\  on the $B_N$-Weyl chamber
\[ \Xi_N = \big\{ x = (x_1,\ldots ,x_N) \in \b R^N: x_1\geq \cdots \geq x_N\geq 0\big\}. \]
The characters of this hypergroup
are given by the $U_N$-means
\begin{equation}\label{Dunklchar} \psi_y(x)= \int_{U_N}
 \mathcal J_\mu\left(\frac{1}{4}y Ux^2 U^{-1}y\right) dU
 =   J_{k(\mu,d)}^B(x, iy), \qquad x,y \in \Xi_N,
 \end{equation}
where
\[k(\mu, d) := \bigl(\mu - (d(N-1)+1)/2,  d/2\bigr)\]
 and elements
from $\Xi_N$ are identif\/ied with diagonal matrices in the natural way.
For details, see Section~4 of~\cite{R2}. We shall now deduce the claimed estimate
for $J_k^B$ from the estimate for the Bessel functions $J_\mu$ on the cone $\Pi_N$ mentioned in the introduction. Indeed, according
to Theorem~3.6 of~\cite{RV2}
there exists a constant $C=C(N,d) >0$ such that for all $\mu > d(2N-1)+1$ and $X\in \Pi_N$,
\begin{equation*}
\big\vert J_\mu(\mu X) - e^{-{\rm tr}\, X}\big\vert  \leq   \frac{C}{\mu} \cdot\min\bigl(1, ({\rm tr}\, X)^2\bigr).
 \end{equation*}
In view of \eqref{Dunklchar}, this leads
to the following estimate for the Dunkl-type Bessel function $J_{k(\mu,d)}^B$:
\[\left\vert J_{k(\mu,d)}^B(2\sqrt{\mu}  x,iy) -
 \int_{U_N} e^{-{\rm tr}(y U x^2U^{-1}y)}dU \right\vert \leq
 \frac{C}{\mu}
\cdot\min \bigl(1, S(x,y) \bigr),
\]
where
\[ S(x,y)  =  \int_{U_N}\big\vert {\rm tr} \big(y Ux^2 U^{-1} y\big)\big\vert ^2 dU  \leq    \int_{U_N} {\rm tr}\big(y^4\big) \cdot  {\rm tr} \big(Ux^4U^{-1}\big) dU  =
 |x^2|^2|y^2|^2 \leq  |x|^4|y|^4.\]
On the other hand,
the Jack polynomials $C_\lambda^\alpha$ with $\alpha = 2/d$  satisfy the product formula
\[ \frac{C_\lambda^\alpha(x^2)C_\lambda^\alpha(y^2)}{C_\lambda^\alpha(\bf{1})}  = \int_{U_N}
C_\lambda^\alpha\big(yUx^2U^{-1}y\big)dU
\qquad \text{for}\quad x,y\in \Xi_N.\]
This follows from a corresponding product formula for the spherical polynomials, see Proposition~5.5 of~\cite{GR}.
 Thus by equation~(\ref{jacksum}) we further obtain, with $\alpha = 2/d$,
\begin{gather} \int_{U_N} e^{-{\rm tr}(y U X^2U^{-1}y)}dU =   \sum_{\lambda\geq 0}
\frac{1}{|\lambda|!} \int_{U_N} C_\lambda^\alpha \big(-y U x^2U^{-1}y\big)dU \nonumber\\
\phantom{\int_{U_N} e^{-{\rm tr}(y U X^2U^{-1}y)}dU }{} =  \sum_{\lambda\geq 0} \frac{1}{|\lambda|!}\frac{C_\lambda^\alpha(-x^2)C_\lambda^\alpha(y^2)}{C_\lambda^\alpha(\bf{1})}
 ={} _0F_0^\alpha \big({-}x^2,y^2\big),\label{harish}
\end{gather}
which implies the assertion.

\begin{remark}
The integral on the left side of formula (\ref{harish}) is of
  Harish-Chandra type.
If $\b F = \b C$, then by Theorem II.5.35 of~\cite{Hel}  it can be written as an alternating sum
\[  \int_{U_N} e^{-{\rm tr}(y U
  x^2U^{-1}y)}dU =\frac{\prod\limits_{j=1}^{N-1}j!}{\pi(x^2)\pi(y^2)}
\sum_{w\in S_N} {\rm sgn}(w) e^{-\langle x^2, wy^2\rangle },\]
where
$ \pi(x) = \prod_{i<j} (x_i-x_j) $ is the fundamental alternating polynomial.
\end{remark}

\section{Proof of  Proposition \ref{proposition11}}\label{sec3}

We now turn to the proof of Proposition \ref{proposition11} which is based on
the power series representa\-tions~(\ref{power-j-a}) and~(\ref{power-j-b}) for the Bessel functions of type A and B. The proof is similar to the
corresponding result for Bessel functions on matrix cones in~\cite{RV2}. We
start with an    observation  about Jack polynomials.

\begin{lemma}\label{lemma20}
For all $x,y\in \b R^N$, $m\in \b N$, and $\alpha>0$,
\[
\sum_{\lambda; |\lambda|=m} \frac{C_\lambda^\alpha(x^2) C_\lambda^\alpha(y^2)}{C_\lambda^\alpha({\bf
    1})} \le  |x|^{2m} |y|^{2m}.
\]
\end{lemma}

\begin{proof} As shown in \cite{KS},
the $C_\lambda^\alpha$ are  nonnegative
linear combinations of monomials.
Therefore,
\[
C_\lambda^\alpha\big(x^2\big)=\sum_{\nu; |\nu|=|\lambda|} c_{\lambda,\nu} x^{2\nu}
\]
with coef\/f\/icients $c_{\lambda,\nu}\ge0, $ and
\[
\frac{C_\lambda^\alpha(x^2)}{C_\lambda^\alpha({\bf 1})} =\sum_{\nu;
  |\nu|=|\lambda|} \tilde c_{\lambda,\nu}  x^{2\nu}
  \]
with suitable $\tilde c_{\lambda,\nu}\ge0$ where $ \sum_{\nu;
  |\nu|=|\lambda|}\tilde c_{\lambda,\nu}=1$.
As $C_\lambda^\alpha(y^2)\ge0 $ and
$ \sum_{ |\lambda|=m}C_\lambda^\alpha\big(y^2\big)=|y|^{2m}$,
we conclude that
\begin{gather*}
\sum_{\lambda; |\lambda|=m} \frac{C_\lambda^\alpha\big(x^2\big) C_\lambda^\alpha\big(y^2\big)}{C_\lambda^\alpha({\bf
    1})}  = \sum_{\lambda,\nu; |\nu|=|\lambda|=m} \tilde c_{\lambda,\nu}
    x^{2\nu}\cdot  C_\lambda^\alpha\big(y^2\big)
\le   |x|^{2m}\cdot \sum_{\lambda,\nu; |\nu|=|\lambda|=m} \tilde c_{\lambda,\nu} C_\lambda^\alpha\big(y^2\big)
\\
\phantom{\sum_{\lambda; |\lambda|=m} \frac{C_\lambda^\alpha(x^2) C_\lambda^\alpha\big(y^2\big)}{C_\lambda^\alpha({\bf
    1})}}{} =  |x|^{2m}\cdot \sum_{\lambda; |\lambda|=m}C_\lambda^\alpha\big(y^2\big)
    =  |x|^{2m} |y|^{2m}
\end{gather*}
as claimed.
\end{proof}

\begin{lemma}\label{lemma2}
Let $\lambda=(\lambda_1,\ldots ,\lambda_N)\ge0$ be a partition, and choose  $k_2\ge0$ and $k_1\ge k_2(N-1)$.
Then for
 $\mu:= k_1 + k_2(N-1)+1/2$, the Pochhammer symbol
$(\mu)_\lambda:= (\mu)_\lambda^{1/k_2}$ satisfies
\[
\left|1-\frac{\mu^{|\lambda|}}{(\mu)_\lambda}\right|\le \frac{1}{3}
2^{N(N-1)(k_2+1)/2}\cdot (1+k_2(N-1))\cdot \frac{|\lambda|^2}{k_1}.
\]
\end{lemma}

\begin{proof}
Consider $  (\mu)_\lambda  =  \prod_{j=1}^N(\mu-k_2(j-1))_{\lambda_j}.$
In this $|\lambda|$-fold product, each factor can be estimated below by
$\mu-k_2(N-1)=k_1+1/2\ge\mu/2$ due to our assumptions. Moreover, precisely
\[ \bigl(0 + 1 + \cdots  + (N-1) \bigr)\lceil k_2\rceil   =   \frac{N(N-1)}{2}\cdot\lceil k_2 \rceil =:r\]
of these factors are smaller than $\mu$. We thus conclude that
\begin{equation*}
(\mu)_\lambda  \ge  \bigl(\mu/2\bigr)^r  \mu^{|\lambda|-r}
 \ge  2^{-N(N-1)(k_2+1)/2} \cdot \mu^{|\lambda|},
\end{equation*}
and thus
\begin{equation}\label{absch-poch}
 \mu^{|\lambda|}/ (\mu)_\lambda  \le   2^{N(N-1)(k_2+1)/2}.
\end{equation}

We next prove by
 induction on the length $|\lambda|$ that
 for $k_1\ge k_2(N-1)$,
\begin{equation}\label{absch}
\left|1-\frac{\mu^{|\lambda|}}{(\mu)_\lambda}\right|\le
\frac{\frac{1}{3}2^{N(N-1)(k_2+1)/2}}{\mu-k_2(N-1)} \cdot (1+k_2(N-1)) \cdot|\lambda|^2.
\end{equation}
As $\mu-k_2(N-1)\ge k_1$, this will imply the lemma.
 In fact, for $k=0,1$, the left
hand side of~(\ref{absch}) is equal to zero, while the right-hand side is nonnegative.
For the induction step, consider a~partition $\lambda$ of length $|\lambda|\ge 2$.
Then there is  a partition $\tilde \lambda$ with $|\tilde \lambda|=|\lambda|-1$ for
which there exists precisely one $j=1,\ldots,N$ with
$\lambda_j=\tilde\lambda_j +1$ while all the other components are equal. Hence, if
we assume the inequality to hold for $ \tilde \lambda$ and use~(\ref{absch-poch}) as well as the abbreviation
$c:= \frac{2}{3}(1+k_2(N-1))$, we obtain
\begin{gather*}
\left| 1-\frac{\mu^{|\lambda|}}{(\mu)_\lambda}\right|=
 \left|1-\frac{\mu^{|\lambda|-1}}{(\mu)_{\tilde\lambda}}
+\frac{\mu^{|\lambda|-1}}{(\mu)_{\tilde\lambda}}
-\frac{\mu^{|\lambda|}}{(\mu)_\lambda}
\right|
\\
\phantom{\left| 1-\frac{\mu^{|\lambda|}}{(\mu)_\lambda}\right|}{} \le
\frac{c}{\mu-k_2(N-1)}\cdot 2^{N(N-1)(k_2+1)/2-1} \cdot (|\lambda|-1)^2\\
\phantom{\left| 1-\frac{\mu^{|\lambda|}}{(\mu)_\lambda}\right|=}{}
 +
\frac{\mu^{|\lambda|-1}}{(\mu)_{\tilde\lambda}}\cdot
\left|1-\frac{\mu}{\mu-k_2(j-1)+\lambda_j-1}\right|
\\
\phantom{\left| 1-\frac{\mu^{|\lambda|}}{(\mu)_\lambda}\right|}{}
\le
\frac{c}{\mu-k_2(N-1)}\cdot 2^{N(N-1)(k_2+1)/2-1} \cdot (|\lambda|-1)^2\\
\phantom{\left| 1-\frac{\mu^{|\lambda|}}{(\mu)_\lambda}\right|=}{}
 +2^{N(N-1)(k_2+1)/2}\cdot\frac{\bigl|-k_2(j-1)+\lambda_j-1\bigr|}{\mu-k_2(j-1)+\lambda_j-1}
\\
\phantom{\left| 1-\frac{\mu^{|\lambda|}}{(\mu)_\lambda}\right|}{}
\le
\frac{2^{N(N-1)(k_2+1)/2-1}}{\mu-k_2(N-1)} \cdot\left(c (|\lambda|-1)^2 +
  2k_2(N-1) + 2|\lambda|-2\right)
\\
\phantom{\left| 1-\frac{\mu^{|\lambda|}}{(\mu)_\lambda}\right|}{}
\le
\frac{2^{N(N-1)(k_2+1)/2-1}}{\mu-k_2(N-1)} \cdot c|\lambda|^2
\end{gather*}
for $|\lambda|\ge2$. Notice that the choice of the constant $c$ is made in
order to ensure that the last inequality holds for $|\lambda|\ge2,$ which
is easily checked by an elementary calculation. This completes the proof.
\end{proof}

We are now ready to prove Proposition \ref{proposition11}.

\begin{proof}[Proof of Proposition \ref{proposition11}]
We use  the power series (\ref{power-j-a}) and  (\ref{power-j-b}) of the
Dunkl--Bessel kernels of type A and B in terms of the
 Jack polynomials $C_\lambda^\alpha$ and the fact that the $C_\lambda^\alpha$
 are homogeneous of degree $|\lambda|$. We thus obtain
\begin{equation*}
 J_{(k_1,k_2)}^B(2\sqrt{\mu}x, iy) - J_{k_2}^A\big({-}x^2,y^2\big)  = \sum_{\lambda\geq 0}
\frac{(-1)^{|\lambda|}}{|\lambda|!}
\left( \frac{\mu^{|\lambda|}}{(\mu)_\lambda}-1\right)
\frac{C_\lambda^\alpha\big(x^2\big) C_\lambda^\alpha\big(y^2\big)}{C_\lambda^\alpha({\bf
    1})}.
\end{equation*}
As
\[
(\mu)_{(1,0,\ldots,0)}=\mu,\qquad (\mu)_{(2,0,\ldots,0)}=\mu(\mu+1),
\qquad \text{and}\qquad
 (\mu)_{(1,1,0,\ldots,0)}=\mu(\mu-k_2),
 \]
the coef\/f\/icients for $|\lambda|\le 1$ are zero, and
 we may write the above expansion as
\[
 J_{(k_1,k_2)}^B(2\sqrt{\mu}x, iy) - J_{k_2}^A\big({-}x^2,y^2\big)  = R_2+R_3
 \]
with
\begin{gather*}R_2=  \frac{1}{2}\left(\frac{\mu^2}{\mu(\mu+1)}-1\right)
\frac{C_{(2,0,\ldots,0)}^\alpha\big(x^2\big)C_{(2,0,\ldots,0)}^\alpha\big(y^2\big)}{C_{(2,0,\ldots,0)}^\alpha({\bf
    1})}
\\
\phantom{R_2=}{}  +  \frac{1}{2}\left(\frac{\mu^2}{\mu(\mu-k_2)}-1\right)
\frac{ C_{(1,1,0,\ldots,0)}^\alpha\big(x^2\big)C_{(1,1,0,\ldots,0)}^\alpha\big(y^2\big)}{C_{(1,1,0,\ldots,0)}^\alpha({\bf
    1})}\end{gather*}
and
\[
R_3=\sum_{m \geq 3} \frac{(-1)^m}{m!} \sum_{|\lambda|= m}
\left( \frac{\mu^{m}}{(\mu)_\lambda}-1\right)
\cdot \frac{C_\lambda^\alpha\big(x^2\big) C_\lambda^\alpha\big(y^2\big)}{C_\lambda^\alpha({\bf
    1})}
.
\]
It now follows from Lemma \ref{lemma20} that under our assumptions on $k_1$ and $ k_2,$
\[
|R_2|\le \frac{C_2}{\mu} \cdot \sum_{|\lambda| =2 } \frac{C_\lambda^\alpha\big(x^2\big) C_\lambda^\alpha\big(y^2\big)}{C_\lambda^\alpha({\bf
    1})}
\le C_2\frac{|x|^4 |y|^4}{\mu}
\]
with some $C_2>0$.
Moreover, Lemmata \ref{lemma2} and  \ref{lemma20} imply that for a suitable constant $C_3$,
\begin{gather*}
|R_3 |
  \le       C_3 \sum_{m\geq 3} \frac{1}{m!} \frac{m^2}{\mu}
\sum_{|\lambda|= m} \frac{C_\lambda^\alpha\big(x^2\big) C_\lambda^\alpha\big(y^2\big)}{C_\lambda^\alpha({\bf
    1})}
\le    C_3 \sum_{m\geq 3} \frac{m^2}{m!\mu}|x|^{2m} |y|^{2m}
\\
\phantom{|R_3 |}{} \le  \frac{2 C_3}{\mu}|x|^4 |y|^4\cdot \sum_{m\geq 3}
 \frac{1}{(m-2)!} |x|^{2m-2} |y|^{2m-2}
 \le  \frac{2 C_3}{\mu}|x|^4 |y|^4\cdot   e^{|x|^2|y|^2 }   .
\end{gather*}
These estimates for $R_2$ and $R_3$ immediately imply the claimed results.
\end{proof}

\pdfbookmark[1]{References}{ref}
\LastPageEnding

\end{document}